\documentclass[a4paper, 12pt]{article}

\usepackage{amsfonts}
\usepackage{amsmath}
\usepackage{graphicx}
\usepackage{multirow}

\newtheorem{theorem}{Theorem}[section]
\newtheorem{lemma}[theorem]{Lemma}
\newtheorem{corollary}[theorem]{Corollary}

\newtheorem{example}[theorem]{Example}

\def\COMPlib{COMP{l$_e\!$ib}}

\title{\LARGE\bf
Hermite matrix in Lagrange basis\\
for scaling static output feedback\\
polynomial matrix inequalities}

\author{Ak{\i}n Deliba{\c s}{\i}$^{1,2,3}$,
Didier Henrion$^{1,2,4}$}

\date{\today}

\begin{document}

\maketitle

\footnotetext[1]{CNRS; LAAS; 7 avenue du colonel Roche, F-31077 Toulouse, France}
\footnotetext[2]{Universit\'e de Toulouse; UPS, INSA, INP, ISAE; LAAS; F-31077 Toulouse, France}
\footnotetext[3]{Department of Electrical Engineering,
Y{\i}ld{\i}z Technical University, Be{\c s}ikta{\c s}, Istanbul, Turkey.}
\footnotetext[4]{Faculty of Electrical Engineering,
Czech Technical University in Prague, Czech Republic.}
\addtocounter{footnote}{4}

\begin{abstract}
Using Hermite's formulation of polynomial stability conditions,
static output feedback (SOF) controller design can be formulated as a
polynomial matrix inequality (PMI), a (generally nonconvex)
nonlinear semidefinite programming problem that can be solved (locally)
with PENNON, an implementation of a penalty method. Typically,
Hermite SOF PMI problems are badly scaled and experiments reveal that
this has a negative impact on the overall performance of the solver.
In this note we recall the algebraic interpretation
of Hermite's quadratic form as a particular B\'ezoutian and we use
results on polynomial interpolation to express the Hermite PMI
in a Lagrange polynomial basis, as an alternative to the conventional
power basis. Numerical experiments on benchmark problem instances
show the substantial improvement brought by the approach, in terms
of problem scaling, number of iterations and convergence behavior
of PENNON.
\end{abstract}

\begin{center}
\end{center}

\section{Introduction}

In 1854 the French mathematician Charles Hermite studied
quadratic forms for counting the number of roots of a polynomial
in the upper half of the complex plane (or, by
a simple rotation, in the left half-plane), more than
two decades before Routh, who was apparently not aware
of Hermite's work, see \cite{jury}. Hurwitz himself used some of Hermite's ideas
to derive in 1895 his celebrated algebraic criterion for
polynomial stability, now called the Routh-Hurwitz criterion and
taught to engineering students in tabular form.

Hermite's criterion can be interpreted
as a symmetric formulation of the Routh-Hurwitz criterion.
This symmetry can be exploited in a semidefinite programming
framework, as shown in \cite{scl99} and \cite{ifip03} in
the context of simultaneous stabilization of linear systems.
Along the same vein, in \cite{cdcecc05} the problem of static
output feedback (SOF) design
was formulated as a polynomial matrix inequality (PMI) problem.
In some cases (e.g. only one input or output available
for feedback) this PMI problem simplifies to a bilinear
matrix inequality (BMI) that can be solved numerically with
PENBMI, a particular instance of PENNON, a general penalty
method for nonlinear and semidefinite programming.
Only convergence to a local optimum is guaranteed, but
experiments reported in \cite{cdcecc05} show that quite
often the approach is viable numerically. In particular,
the SOF PMI formulation involves only controller parameters,
and does not introduce (a typically large number of)
Lyapunov variables.

Our motivation in this paper is to contribute along the lines
initiated in \cite{cdcecc05} and to study the impact of
SOF PMI problem formulation on the behavior of PENNON,
in particular w.r.t. data scaling and number of iterations.
The Hermite matrix depends quadratically on coefficients
of the characteristic polynomial, in turn depending
polynomially on the controller parameters. As a result,
coefficients of a given Hermite matrix typically differ by
several orders of magnitude, and experiments reveal
that this poor data scaling significantly impacts
on the performance of PENNON.

In this paper we use an alternative formulation of the
Hermite matrix, using a Lagrange polynomial basis instead
of the standard power basis. We build on previous work
from the computer algebra and real algebraic geometry
communities, recalling the interpretation
of Hermite's quadratic form as a particular B\'ezoutian,
the resultant of two polynomials, see \cite{issac08}
and references therein.
This interpretation provides a natural choice for
the nodes of the Lagrange basis. The construction
of the Hermite matrix in this basis is carried out
efficiently by interpolation, overcoming difficulties
inherent to Vandermonde matrices, as suggested in \cite{shakoori}
for general B\'ezout matrices.

In addition to digesting and tailoring to our needs
results from computational algebraic geometry, another
contribution of our paper is to extend slightly the
characterization of \cite{shakoori} to Hermitian
forms with complex and repeated interpolation nodes.
In particular, in our SOF design application framework,
these nodes are roots of either imaginary or real part of a target characteristic polynomial
featuring spectral properties desirable for the
closed-loop system. This target polynomial is the
main tuning parameter of our approach, and we provide
numerical evidence that a suitably choice of target
polynomial, compatible with achievable closed-loop
dynamics, results in a significant improvement
of SOF PMI problem scaling, with positive effects
on the overall behavior (convergence, number of
outer and inner iterations, linesearch steps) of PENNON.
Furthermore, some of the problems that were not solvable
in the power basis, see \cite{cdcecc05}, can now be solved
in the Lagrange basis. These improvements are illustrated
on numerical examples extracted from the publicly available
benchmark collection \mbox{\COMPlib}, see \cite{compleib}.

\section{PMI formulation of SOF design problem}

We briefly recall the polynomial matrix inequality (PMI)
formulation of static output feedback (SOF) design
problem proposed in \cite{cdcecc05}.

Consider the linear system
\[
\begin{array}{rcl}
\dot{x} & = & Ax+Bu \\
y & = & Cx
\end{array}
\]
of order $n$ with $m$ inputs and $p$ outputs, that
we want to stabilize by static output feedback (SOF)
\[
u = Ky.
\]
In other words, given matrices $A \in {\mathbb R}^{n \times n}$,
$B \in {\mathbb R}^{n \times m}$, $C \in {\mathbb R}^{p \times n}$,
we want to find matrix $K \in {\mathbb R}^{m \times p}$ such
that the eigenvalues of closed-loop matrix $A+BKC$ all belong
to the left half of the complex plane.

Let $k \in {\mathbb R}^{mp}$ be the vector obtained by stacking
the columns of matrix $K$. Define
\begin{equation}\label{charpol}
q(s,k) = \mathrm{det}\:(sI-A-BKC) = \sum_{i=0}^n q_i(k)s^i
\end{equation}
as the characteristic polynomial of matrix $A+BKC$.
Coefficients of increasing powers of indeterminate $s$ in
polynomial $q(s,k)$ are multivariate polynomials in $k$, i.e.
\begin{equation}\label{coefs}
q_i(k) = \sum_{\alpha} {q_i}_{\alpha} k^{\alpha}
\end{equation}
where $\alpha \in {\mathbb N}^{mp}$ describes all monomial powers.

The Routh-Hurwitz criterion for stability of polynomials has a symmetric version called
the Hermite criterion. A polynomial is stable if and only if its Hermite matrix, quadratic
in the polynomial coefficients, is positive definite. Algebraically, the Hermite matrix
can be defined via the B\'ezoutian, a symmetric form of the resultant.

Let $a(u)$, $b(u)$ be two polynomials of degree $n$ of the indeterminate $u$.
Define the bivariate quadratic form
\[
\frac{a(u)b(v)-a(v)b(u)}{u-v}=\sum_{i=1}^n\sum_{j=1}^n b_{ij}u^{i-1}v^{j-1}.
\]
The $n$-by-$n$ matrix with entries $b_{ij}$ is the B\'ezoutian matrix,
whose determinant is the resultant of $a$ and $b$,
obtained by eliminating variable $u$ from the system of equations
$a(u)=b(u)=0$.

The Hermite matrix in power basis of $q(s,k)$, denoted by $H^P(k)$,
is defined as the B\'ezoutian matrix of the real and imaginary parts of $q(ju,k)$:
\[
\begin{array}{ccc}
a(u,k) &=&\mathrm{Im}\,q(jw,k)\\
b(u,k)&=&\mathrm{Re}\,q(jw,k).
\end{array}
\]
The roots of polynomial $q(s,k)$ belongs to the left half-plane
if and only if
\[
H^P(k) = \sum_{i=0}^n \sum_{j=0}^n q_i(k)q_j(k) H^P_{ij} \succ 0.
\]
The above relation is a matrix inequality depending polynomially on parameters $k$.
Therefore, finding $k$ amounts to solving a polynomial matrix
inequality (PMI) problem.

\begin{example}\label{nn6}\rm
As an illustrative example, consider problem $\tt NN6$ in \cite{compleib}.
The closed-loop characteristic polynomial is (to 8 significant
digits):
\begin{align*}
q(s,k)\:\:=\:\: & s^9+23.300000s^8+(4007.6500-14.688300k_2+14.685000k_4)s^7\\
&+(91133.935-14.685000k_1+14.688300k_3+15.132810k_4)s^6\\
&+(1149834.9-57334.489k_2+15.132810k_3+36171.693k_4)s^5\\
&+(20216420-57334.489k_1+36171.693k_3+35714.763k_4)s^4\\
&+(49276365-12660338k_2+35714.763k_3+3174671.8k_4)s^3\\
&+(-1562.6281\cdot 10^5-12660338k_1-3174671.8k_3+3133948.9k_4)s^2\\
&+(-4315.5562\cdot 10^5+95113415k_2+3133948.9k_3)s\\
&+95113415k_1\\
\end{align*}
with SOF gain $K=[k_1\,k_2\,k_3\,k_4]$. The $9$-by-$9$ Hermite
matrix of this polynomial cannot be displayed entirely for
space reasons, so we choose two representative entries:
\begin{align*} H^P_{3,3}(k)\:\:=\:\:&10244466\cdot 10^8-53923375\cdot 10^7k_1+55487273\cdot 10^6k_2\\
&+10310826\cdot 10^7k_3-32624061\cdot 10^7k_4+16028416\cdot 10^7k_1k_2\\
&-27103829\cdot 10^4k_1k_3-36752006\cdot 10^6k_1k_4\\
&-43632833\cdot 10^6k_2k_3-43073807\cdot 10^6k_2k_4\\
&+22414163k_3^2+10078541\cdot 10^6k_3k_4+99492593\cdot 10^5k_4^2
\end{align*}
and
\[
H^P_{9,9}(k) = 23.300000.
\]
We observe that this Hermite matrix is ill-scaled, in the sense
that the coefficients of its entries (multivariate polynomials in $k_i$)
differ by several orders of magnitude. This representation is not suitable
for a matrix inequality solver.
\end{example}

\section{A simple scaling strategy}

A possible remedy to address the poor scaling properties
of the Hermite matrix is to scale the frequency variable $s$, that is,
to substitute $\rho s$ for $s$ in the characteristic polynomial $q(s,k)$,
for a suitable positive scaling $\rho$. Finding the optimal value
of $\rho$ (e.g. in terms of relative scaling of the coefficients
of the Hermite matrix) may be formulated as an optimization problem,
but numerical experiments indicate that nearly optimal results are
achieved when following the basic strategy consisting
of choosing $\rho$ such that the constant and highest
power polynomial coefficients are both equal to one.
For example, this idea was implemented by Huibert Kwakernaak
in the {\tt scale} function
of the Polynomial Toolbox for Matlab, see \cite{polyx}.

\begin{example}\label{ac4}\rm
Consider the simple example $\tt AC4$ in \cite{compleib}.
The open-loop characteristic polynomial is
\begin{align*}
q(s,0)&=\mathrm{det}(sI-A)\\
&=s^4+150.92600s^3+130.03210s^2-1330.6306s-66.837750
\end{align*}
with Hermite matrix in power basis
\[
H^P=\left [
\begin{array}{cccc}
88936.354&0&10087.554&0\\
0&-162937.14&0&1330.631\\
10087.554&0&20955.855&0\\
0&1330.6306&0&150.92600
\end{array}
\right ].
\]
To measure quantitatively the scaling of a matrix $X$,
we may use its condition number. If the matrix is poorly
scaled, then its condition number is large. Minimizing
the condition number therefore improves the scaling.
For the above matrix, its condition number (in the Frobenius
norm), defined as $\Vert H^P \Vert_F \Vert (H^P)^{-1} \Vert_F$,
is equal to $1158.2$. If we choose $\rho=\sqrt[4]{1.0000/66.840}=3.5000\cdot 10^{-1}$,
the scaled characteristic polynomial has unit constant
and highest coefficient, and the resulting scaled Hermite
matrix reads
\[
S H^P S=\left [
\begin{array}{cccc}
163.48864&0&151.37636&0\\
0&-2445.0754&0&163.00225\\
151.37636&0&2567.0923&0\\
0&163.00225&0&150.92600
\end{array}
\right ]
\]
with
\[
S=\mathrm{diag} \left(\rho^3,  \rho^2, \rho^1, 1 \right).
\]
The Frobenius condition number of $S H^P S$ is equal to $32.096$.
\end{example}

Whereas this simple scaling strategy with one degree of freedom
may prove useful for small-degree polynomials and small-size
Hermite matrices, a more sophisticated approach is required
for larger instances.

\section{Hermite matrix in Lagrange basis}\label{hermite_cr_lag}

In this section we show how the Hermite matrix can be scaled
by an appropriate choice of polynomial basis. Moreover, this basis
allows for a straightforward entrywise construction of
the Hermite matrix.

\subsection{Distinct interpolation points}

Consider $n$ distinct interpolation
points $u_i \in {\mathbb C}$, $i=1,\ldots,n$, and define the $j$-th Lagrange polynomial
\[
l_j(u) = \prod_{i=1,i \neq j}^{n} \frac{u-u_i}{u_j-u_i}
\]
which is such that $l_j(u_j)=1$ and $l_j(u_i)=0$ if $i\neq j$.
In matrix form we can write
\begin{equation}\label{vandermonde}
\left[\begin{array}{c}1\\u\\u^2\\\vdots\\u^{n-1}\end{array}\right] =
\left[\begin{array}{cccc}1&1&\cdots&1\\u_1&u_2&\cdots&u_n\\u_1^2&u_2^2&\cdots&u_n^2\\
\vdots&&&\vdots\\u_1^{n-1}&u_1^{n-1}&\cdots&u_n^{n-1}\end{array}\right]
\left[\begin{array}{c}l_1(u)\\l_2(u)\\l_3(u)\\\vdots\\l_n(u)\end{array}\right] =
V_u l(u)
\end{equation}
where $V_u$ is a Vandermonde matrix.
Given a univariate polynomial $q(s)$ with real coefficients, define
\begin{equation}\label{realimag}
\begin{array}{ccc}
a(u) &=&\mathrm{Im}\,q(jw)\\
b(u)&=&\mathrm{Re}\,q(jw)
\end{array}
\end{equation}
as its imaginary and real parts on the imaginary axis, respectively.
In the following, the star denotes transpose conjugation and the
prime denotes differentiation, i.e.
\[
a'(u) = \frac{da(u)}{du}.
\]

\begin{theorem}\label{lagrange}
When the interpolation points are distinct (i.e. $u_i \neq u_j$, $i \neq j$, $i,j=1,\ldots,n$),
the Hermite matrix of $q(s)$ in Lagrange basis, denoted
by $H^L$, is given entrywise by
\[
H^L_{i,j} :=\left\{
\begin{array}{rl}
\displaystyle\frac{a(u_i^*)b(u_j)-a(u_j)b(u_i^*)}{u_i^*-u_j} & \mathrm{if}\:u_i^* \neq u_j,\\
a'(u_i^*)b(u_j)-a(u_j)b'(u_i^*) & \mathrm{otherwise,}
\end{array}
\right.
\]
for all $i,j=1,\ldots,n$.
\end{theorem}

{\bf Proof}
Let us express the B\'ezoutian of $a$ and $b$ as a bivariate quadratic form
\[
\frac{a(u)b(v)-a(v)b(u)}{u-v}=\left[\begin{array}{c}1 \\ v \\ \vdots\\ v^{n-1}\end{array}\right ]^*
H^P \left [ \begin{array}{c} 1 \\ u \\ \vdots \\ u^{n-1} \end{array} \right ]
\]
where $H^P$ is the Hermite matrix of $q$
in the power basis. Recalling relation (\ref{vandermonde}), the B\'ezoutian becomes
\[
\frac{a(u)b(v)-a(v)b(u)}{u-v}=l(v)^*V_v^*H^PV_ul(u) = l(v)^*H^Ll(u)
\]
so that the Hermite matrix of $q$ in the Lagrange basis can be expressed as
\[
H^L = V_v^*H^PV_u.
\]
By evaluation at $n$ distinct interpolation points $u_i$ and $v_j$,
$H^L$ is given entrywise by
\begin{equation}\label{entry}
H^L_{i,j} = \frac{a(u_i^*)b(v_j)-a(v_j)b(u_i^*)}{u_i^*-v_j}.
\end{equation}
Now let $u_i^*\rightarrow v_j$ for all $i,j=1,\ldots,n$.
After adding and subtracting $a(v_j)b(v_j)$ to the numerator of (\ref{entry}), we find
\[
H^L_{i,j}=a'(u_i^*)b(u_j)-a(u_j)b'(u_i^*),
\]
using a limiting argument.
$\Box$

\subsection{Repeated interpolation points}

Let us define the bivariate polynomials
\[
c_{i,j}(u,v):=\frac{\partial^{i+j-2}}{\partial u^{i-1}\partial v^{j-1}}\left (\frac{a(u)b(v)-a(v)b(u)}{u-v}\right)
\]
for all $i,j=1,\ldots,n$ and denote by
\[
a^{(k)}(u) = \frac{d^ka(u)}{du^k}
\]
the $k$-th derivative of univariate polynomial $a(u)$.

\begin{lemma}\label{rep_Lagrange}
When the interpolation points are all equal (i.e. $u_i=u_j$ for all $i,j=1,\ldots,n$),
the Hermite matrix of $q(s)$ in Lagrange basis is given entrywise by
\[
H^L_{i,j}:=\left\{
\begin{array}{lr}
\displaystyle{\frac{c_{i,j}(u_i^*,u_j)}{(i-1)!(j-1)!}}& \mathrm{if}\:u_i^* \neq u_j,\vspace{5mm}\\ 
\displaystyle\sum_{k=0}^{i-1}\frac{a^{(j+k)}(u_i^*)b^{(i-k-1)}(u_j)-a^{(i-k-1)}(u_j)b^{(j+k)}(u_i^*)}{(j+k)!(i-k-1)!}& \mathrm{otherwise,}
\end{array}
\right.
\]
for all $i,j=1,\ldots,n$.
\end{lemma}
{\bf Proof}
The proof of this result follows along the same lines
as the proof of Theorem \ref{lagrange}, with additional
notational difficulties due to higher-order differentations.
$\Box$

\begin{example}\label{rep_HL}\rm
Let us choose $n=3$ equal interpolation points ($u_1=u_2=u_3=x\in{\mathbb R}$).
According to Lemma \ref{rep_Lagrange}, $H^L$ has the following entries:
\[
\begin{array}{rcl}
H_{11}^L &=& \frac{a'(x)b(x)-a(x)b'(x)}{1!} \\
H_{12}^L &=& \frac{a^{(2)}(x)b(x)-a(x)b^{(2)}(x)}{2!}\\
H_{13}^L &=& \frac{a^{(3)}(x)b(x)-a(x)b^{(3)}(x)}{3!}\\
H_{22}^L &=& \frac{a^{(2)}(x)b'(x)-a'(x)b^{(2)}(x)}{2!}+\frac{a^{(3)}(x)b(x)-a(x)b^{(3)}(x)}{3!}\\
H_{23}^L &=& \frac{a^{(3)}(x)b'(x)-a'(x)b^{(3)}(x)}{3!}+\frac{a^{(4)}(x)b(x)-a(x)b^{(4)}(x)}{4!}\\
H_{33}^L &=& \frac{a^{(3)}(x)b^{(2)}(x)-a^{(2)}(x)b^{(3)}(x)}{3!2!}+\frac{a^{(4)}(x)b'(x)-a'(x)b^{(4)}(x)}{4!}+\frac{a^{(5)}(x)b(x)-a(x)b^{(5)}(x)}{5!}.
\end{array}
\]
\end{example}

Based on Theorem \ref{lagrange} and Lemma \ref{rep_Lagrange},
we leave it to the reader to derive entrywise expressions for the Lagrange basis
Hermite matrix in the general case when only some interpolation points are repeated.

In the remainder of the paper we will assume for notational simplicity that the
interpolation points are all distinct.

\subsection{Scaling}

\begin{corollary}\label{diag}
Let the interpolation points be (distinct) roots of either $a(u)$ or $b(u)$,
as defined in (\ref{realimag}).
Then the Hermite matrix of $q(s)$ in Lagrange basis is block diagonal, with $2\times 2$ blocks corresponding to pairs of complex conjugate points
and $1\times 1$ blocks corresponding to real points.
\end{corollary}

{\bf Proof}
From Theorem \ref{lagrange}, all the off-diagonal entries of $H^L$ are given by
\begin{equation}\label{off_diag_herm}
\displaystyle\frac{a(u_i^*)b(u_j)-a(u_j)b(u_i^*)}{u_i^*-u_j}
\end{equation}
when interpolation points $u_i$ and $u_j$ are not complex conjugate.
Both terms $a(u_i^*)b(u_j)$ and $a(u_j)b(u_i^*)$ are equal to zero in (\ref{off_diag_herm}) since the interpolation points are the roots of either $a(u)$ or $b(u)$. The diagonal entries are
$a'(u_i^*)b(u_j)-a(u_j)b'(u_i^*)$ since it is assumed that interpolation points are distinct.
Therefore this part of $H^L$ is $1\times 1$ block-diagonal.

When interpolation points $u_i$ and $u_j$ are complex conjugate, there is only one non-zero entry
$(i,j)$ which is equal to $a'(u_i^*)b(u_j)-a(u_j)b'(u_i^*)$ and located in the off-diagonal entry, according to pairness. The diagonal entries of this case are equal to zero by virtue of equation (\ref{off_diag_herm}). Therefore this part of $H^L$ is $2\times 2$ block-diagonal.$\Box$

From Corollary \ref{diag} it follows that we can easily find a block-diagonal scaling matrix
$S$ such that the scaled Lagrange Hermite matrix
\[
H^S = S H^L S
\]
has smaller condition number. Nonzero entries of $S$ are given by
\[
S_{i,j}:=\left|\sqrt{H^L_{i,j}}\right|^{-1}
\]
whenever $H^L_{i,j}$ is a nonzero entry $(i,j)$ of $H^L$.

\begin{example}\label{nn5}\rm
As an illustrative example, consider problem $\tt NN5$ in \cite{compleib}.
The open-loop characteristic polynomial is
\[
\begin{array}{rcl}
q(s) & = & s^7+10.171000s^6+96.515330s^5+458.42510s^4\\
& & +2249.4849s^3+1.2196400s^2-448.72180s+6.3000000.
\end{array}
\]
The Hermite matrix in power basis has the following entries:
\[
\begin{array}{rclrcl}
	H^P_{1,1}&=&-2826.9473&H^P_{1,3}&=&-14171.755\\
	H^P_{1,5}&=&608.04658&H^P_{1,7}&=&-6.3000000\\
        H^P_{2,2}&=&-14719.034&H^P_{2,4}&=&206313.38\\
	H^P_{2,6}&=&-4570.2494&H^P_{3,3}&=&209056.94\\
	H^P_{3,5}&=&-4687.9634&H^P_{3,7}&=&1.2196400\\
        H^P_{4,4}&=&1026532.4&H^P_{4,6}&=&-22878.291\\
	H^P_{5,5}&=&21366.759&H^P_{5,7}&=&-458.42510\\
	H^P_{6,6}&=&523.23232&H^P_{7,7}&=&10.171000,\\
\end{array}
\]
remaining nonzero entries being deduced by symmetry. Apparently, this matrix
is ill-scaled. Choosing interpolation points $u_i$ as roots of $a(u)$,
the imaginary part of $q(s)$ along the imaginary axis, we
use Theorem \ref{lagrange} to build the Hermite matrix in
Lagrange basis:
\[
\begin{array}{rcl}
H^L & = & \mathrm{diag}\:(
-2826.9473,\: 41032866\cdot10^3,\: 44286011\cdot10^2,\:
41032866\cdot10^3,\\
& & \quad 44286011\cdot10^2,\:\left[\begin{array}{cc}
0 & 22222.878\\ 22222.878 & 0\end{array}\right]).
\end{array}
\]
This matrix
is still ill-scaled (with Frobenius condition number equal to $2.0983\cdot 10^7$),
but it is almost diagonal.
Using an elementary
diagonal scaling matrix $S$, we obtain
\[
H^S=S H^L S=\mathrm{diag}\:(
-1,\:1,\:1,\:1,\:1,\:\left[\begin{array}{cc}0&1\\1&0\end{array}\right])
\]
which is a well-scaled representation of the Hermite matrix,
with Frobenius condition number equal to $7$.
\end{example}

\subsection{Target polynomial}

In our control application, let us introduce our main tuning tool which we call target polynomial,
denoted by $q(s)$.
The target polynomial provides the interpolation points required to build well-scaled Hermite matrix
in the SOF problem. These points are defined as in Corollary \ref{diag} as the roots of
either the real or imaginary part of $q(s)$ when evaluated along the imaginary axis.

In the context of SOF design, the target polynomial may be either choosen as
\begin{itemize}
\item a valid closed-loop characteristic polynomial (\ref{charpol}) for a specific value of $k$, or
\item a polynomial with desired pole distribution for the closed-loop system.
\end{itemize}

Furthermore, we invoke a continuity argument to observe that
the condition and/or scaling of the Hermite matrix does not change
abruptly in a neighborhood of a given target polynomial.

\begin{example}\label{nn6b}\rm
Consider again Example \ref{nn6} and let the target polynomial be an achievable
closed-loop characteristic polynomial $q(s,k)=\mathrm{det}(sI-A-BKC)$,
where
\[
K=[-4.3264\cdot 10^{-1},\:\:-1.6656,\:\:1.2537\cdot 10^{-1},\:\:2.8772\cdot 10^{-1}]
\]
is a random feedback gain.
The roots of the imaginary part of $q(s,k)$ are chosen as interpolation points
\[
u=\left(0,\:\pm60.847,\:\pm16.007,\:\pm9.2218,\pm2.7034i\right).
\]
Here are two representative entries
of the resulting Lagrange basis Hermite matrix:
\begin{align*}
H^S_{3,3}(k)\:\:=\:\:&9.4439251\cdot 10^{-1}+1.9763715\cdot 10^{-4}k_1-8.9049916\cdot 10^{-4}k_2\\
&-8.6909277\cdot 10^{-3}k_3+1.9212126\cdot 10^{-1}k_4\\
&+3.8300306\cdot 10^{-9}k_1k_2-1.0276186\cdot 10^{-8}k_1k_3\\
&+3.3905595\cdot 10^{-5}k_1k_4-3.4222179\cdot 10^{-5}k_2k_3\\
&-3.8046300\cdot 10^{-5}k_2k_4+2.7420115\cdot 10^{-9}k_3^2\\
&+6.5442491\cdot 10^{-6}k_3k_4+1.015195648\cdot 10^{-5}k_4^2\\
\end{align*}
and
\[
H^S_{1,1}(k)=-1.6918611k_1+3.7288052\cdot 10^{-1}k_1k_2+1.2286264\cdot 10^{-2}k_1k_3.
\]
Comparing with the entries of the power basis Hermite matrix $H^P(k)$
given in Example \ref{nn6}, we observe
a significant improvement in terms of coefficient scaling.
\end{example}

\section{Numerical examples}

In this section, we present the benefits of Lagrange basis against power basis
when solving SOF PMI problems found in the database \mbox{\COMPlib}, see \cite{compleib}.
Even though Michal Ko{\v c}vara and Michael Stingl informed us that an AMPL
interface to PENNON is now available to solve PMI problems, in this paper
for simplicity we consider only BMIs (i.e. quadratic PMIs) and the
PENBMI solver (a particular instance of PENNON focusing on BMIs)
under the YALMIP modeling interface, see \cite{yalmip}.
The numerical examples are processed with YALMIP R20070523 and PENBMI 2.1
under Matlab R2007a running on a Pentium D 3.4GHz system with 1GB ram.
We set the PENBMI penalty parameter {\tt P0} by default to
$0.001$ (note that this is not the default YALMIP setting).

As in \cite{cdcecc05}, the optimization problem to be solved is
\[
\begin{array}{ll}
\min_{k,\lambda} & \mu \|k\| - \lambda \\
\mathrm{s.t.} & H(k) \succeq \lambda I
\end{array}
\]
where $H(k)$ is the Hermite matrix in power or Lagrange basis,
$\mu > 0$ is a parameter and $\|.\|$ is the Euclidean norm.
Parameter $\mu$ allows to trade off
between feasibility of the BMI and a moderate norm of the feedback
gain, which is generally desirable in practice, to avoid
large feedback signals. This adjustment is necessary in many examples.
Indeed, the smallest values of $\|k\|$ are typically located at the boundary
of the feasibility set, so the resulting closed-loop system
is fragile and a small perturbation on system parameters may be
destabilizing.

PENBMI is a local optimization solver. Therefore, the choice of
initial guess $k_0$, $\lambda_0$ is critical. In most of the examples
we choose the origin as the initial point. However this is not always
an appropriate choice, as illustrated below. In addition to this,
PENBMI does not directly handle complex numbers (unless the real
and imaginary parts are split off, resulting in a real coefficient
problem of double size), so we restrict the interpolation points
to be real numbers.

As a result of the root interlacing property, the roots of real
and imaginary parts of a stable polynomial are real (and interlacing).
Owing to this fact, if we choose a stable target polynomial $q(s)$
the resulting interpolation points are necessarily real.

\begin{example}\label{ac4b}\rm
Consider again problem $\tt AC4$, with characteristic polynomial
\begin{align*}
q(s,k)\:\:=\:\:&s^4+150.92600s^3+(130.03210-18.135000k_1-19612.500k_2)s^2\\
&-(1330.6306+19613.407k_1+18322.789k_2)s-(66.837750+980.62500k_1+867.10818k_2)
\end{align*}
and power basis Hermite matrix with entries
\begin{align*}
H^P_{1,1}\:\:=\:\:&88936.354+2615765.6k_1+2378454.6k_2\\
&+19233397k_1^2+34974730k_1k_2+15887840k_2^2\\
H^P_{1,3}\:\:=\:\:&10087.554+148001.81k_1+130869.17k_2\\
H^P_{2,2}\:\:=\:\:&-162937.14-2378239.7k_1+23845311k_2\\
&+355689.13k_1^2+38500022\cdot 10^{1}k_1k_2+35935569\cdot 10^{1}k_2^2\\
H^P_{2,4}\:\:=\:\:&1330.6306+19613.407k_1+18322.789k_2\\
H^P_{3,3}\:\:=\:\:&20955.855+16876.364k_1-2941713.4k_2\\
H^P_{4,4}\:\:=\:\:&150.92600.\\
\end{align*}
Open-loop poles of the system are $(2.5792\:,-5.0000\cdot10^{-2}\:,-3.4552,\:-150.00)$.
If we define our target polynomial roots as $(-5.0000\cdot10^{-2}\:,-5.0000\cdot10^{-2}\:,-3.4552\:,-150.00)$,
keeping the stable open-loop poles and shifting the unstable open-loop pole
to the left of the imaginary axis, our 4 interpolation points (roots of
the real part of the target polynomial)
are $u=(\pm 23.100\:,\pm4.9276\cdot10^{-2})$ and the resulting
Lagrange basis Hermite matrix has entries
\begin{align*}
H^S_{1,1}\:\:=\:\:&6.3432594\cdot 10^{-1}+3.1878941\cdot10^{-1}k_1-17.462079k_2\\
&+4.4822907k_1^2+4.4060140k_1k_2+4.1121739k_2^2\\
H^S_{1,2}\:\:=\:\:&-3.7795293\cdot 10^{-1}-1.0581354k_1-18.455273k_2\\
&-3.6574685\cdot 10^{-3}k_1^2-4.4045142k_1k_2-4.1114926k_2^2\\
H^S_{1,3}\:\:=\:\:&2.4459288+36.285639k_1+42.619220k_2\\
&+7.8459729k_1^2+189.06139k_1k_2+169.77324k_2^2\\
H^S_{1,4}\:\:=\:\:&1.9481147+28.929191k_1+12.037605k_2\\
&+7.5224565k_1^2-161.11487k_1k_2-157.07808k_2^2\\
H^S_{2,2}\:\:=\:\:&6.3432594\cdot10^{-1}+3.1878941\cdot10^{-1}k_1-17.462079k_2\\
&+4.4822907\cdot10^{-3}k_1^2+4.4060140k_1k_2+4.1121739k_2^2\\
\end{align*}
\begin{align*}
H^S_{2,3}\:\:=\:\:&1.9481074+28.929083k_1+12.037568k_2\\
&+7.5224134k_1^2-161.11415k_1k_2-157.07737k_2^2\\
H^S_{2,4}\:\:=\:\:&2.4459288+36.285639k_1+42.619220k_2\\
&+7.8459729k_1^2+189.06139k_1k_2+169.77324k_2^2\\
H^S_{3,3}\:\:=\:\:&659.47243+19434.408k_1+18140.083k_2\\
&+143181.92k_1^2+267314.59k_1k_2+124766.18k_2^2\\
H^S_{3,4}\:\:=\:\:&665.36241+19520.378k_1+17279.067k_2\\
&+143169.06k_1^2+253396.76k_1k_2+111775.41k_2^2\\
H^S_{4,4}\:\:=\:\:&659.47243+19434.408k_1+18140.083k_2\\
&+143181.92k_1^2+267314.59k_1k_2+124766.18k_2^2.
\end{align*}
Choosing the power basis representation with the orgin as initial point and
trade-off parameter $\mu=10^{-5}$, PENBMI stops by a linesearch failure and YALMIP displays
a warning. However, we obtain a feasible solution $\lambda=150.88$ and $K=[1.4181,\:\:-1.6809]$.
This computation requires $43$ outer iterations, $433$ inner iterations and $825$ linesearch steps.
On the other hand, in the Lagrange basis representation, the problem
was solved with no error or warning, yielding $\lambda=9.8287\cdot 10^{-1}$, $K=[-5.0902\cdot 10^{-2},\:\:-2.0985\cdot 10^{-2}]$
with
$17$ outer iterations, $100$ inner iterations and $159$ linesearch steps.

We notice however that using the same trade-off parameter $\mu$ for both representations is not fair
since $H^P$ and $H^S$ have significantly different scalings. If we choose $\mu=0.1$ for
the power basis representation, no problem is detected during the process and we obtain
$\lambda=150.87$, $K=[8.0929\cdot 10^{-2},\:\:-1.6953\cdot 10^{-1}]$
after $26$ outer iterations, $188$ inner iterations and
$238$ linesearch steps. So it seems that the Lagrange basis representation becomes
relevant mainly for high degree systems. This is confirmed by the experiments below.
\end{example}

Consider the $\tt AC7$, $\tt AC17$, $\tt REA3$, $\tt UWV$, $\tt NN5$, $\tt NN1$ and $\tt HE1$
SOF BMI problems of \COMPlib. In Table \ref{tab:result} we report comparative results
for the power and Lagrange basis representations. As in Example \ref{ac4},
the main strategy to choose the target polynomials (and hence the interpolation points)
is to mirror the open-loop stable roots, and to shift the open-loop roots to, say $-5.0000\cdot10^{-1}$
(any other small negative value may be suitable).
We see that the behavior indicators of PENBMI are significantly better in the
Lagrange basis, and the improvement is more dramatic for larger degree
examples. More specifically:
\begin{itemize}
\item for small degree systems like $\tt AC17$
there is only a minor improvement;
\item at the first attempt to solve the $\tt REA3$ example
strict feasibility was not achieved in the power basis, since $\lambda$ is almost zero.
Therefore it was necessary to tune the $\mu$ parameter. Results of the second attempt
show that the BMI problem was solved and the Lagrange basis computation was slightly
less expensive than the power basis computation;
\item the underwater vehicle example $\tt UWV$ has two inputs and two outputs.
However, because of cancellation of higher degree terms in the characteristic
polynomial, the degree of the Hermite matrix is equal to $2$ and we can use PENBMI
on this problem;
\item on open-loop stable systems such as $\tt UWV$ or $\tt AC17$, the
improvement brought by the Lagrange basis is less significant. Since the
main purpose of our optimization problem is to minimize the norm of control gain,
we observe that the Lagrangian basis is still slightly better than the power basis;
\item PENBMI is unable to reach a feasible point for examples $\tt NN5$, $\tt NN1$ and $\tt HE1$,
when we choose the origin as the initial point. Indeed, local optimization techniques seek an
optimal point inside the feasible set in a neighborhood of the initial point. Therefore, achievement of the solver may be very sensitive to the initial point. When the initial point is defined heuristically or
randomly, the improvement is significant for system $\tt NN5$ in Lagrange basis.
However, there is no improvement over $\tt NN1$ and $\tt HE1$, when we use this simple strategy to define the target polynomial.
\end{itemize}

\begin{table}[h]
\setlength{\tabcolsep}{0.8pt}
\caption {PENBMI performance on SOF BMI problems}\label{tab:result}
\begin{center}
\small
\begin{tabular}{|c|c|c|c|c|c|c|c|c|}
    \hline
	\multirow{2}{*}{system} & \multirow{2}{*}{basis} & \multirow{2}{*}{$\mu$} & \multirow{2}{*}{$K_0$} & out.& inn.& lin.& \multirow{2}{*}{$K$} & \multirow{2}{*}{$\lambda$}\\ 
	&&&&iter.&iter.&steps&&\\
	\hline
	\multirow{2}{*}{$\begin{array}{c} \tt AC7 \\ n=9\end{array}$ }& pow. & $1$ &$[0\,0]$ & $27$ & $148$ & $167$ & $[1.1205\:\:\:-3.0946\cdot 10^{-1}]$ & $51.640$\\
\cline{2-9}
    & Lag. & $10^{-5}$ & $[0\,0]$ & $15$ & $51$ & $67$ & $[5.7336\:\:\:3.9995]$ & $3.6356\cdot 10^{-1}$\\
	\hline
	\multirow{2}{*}{$\begin{array}{c} \tt AC17 \\ n=4\end{array}$ }& pow. & $1$ &$[0\,0]$ & $14$ & $65$ & $173$ &  $[1.6619\cdot 10^{-1}\:\:\:8.5782\cdot 10^{-1}]$ & $5.8306$\\
\cline{2-9}
    & Lag. & $1$ & $[0\,0]$ & $16$ & $36$ & $57$ & $[-1.0855\cdot 10^{-2}\:\:\:1.5128\cdot 10^{-1}]$ & $1.0459$\\
	\hline
	\multirow{3}{*}{$\begin{array}{c} \tt REA3 \\ n=12\end{array}$ }& \multirow{2}{*}{pow.} & $1$ &$[0\,0\,0]$ & $21$ & $28$ & $28$ &  $\left [ \begin{array}{c}0\:\:-1.0435\cdot 10^{-5}\:\:-2.2281\cdot 10^{-4}\end{array}\right ]$ & $8.4187\cdot 10^{-13}$\\
\cline{3-9}
    & & $10^{-5}$ & $[0\,0\,0]$ & $46$ & $458$ & $2460$ & $\left [ \begin{array}{c}0\:\:-43711\:\:-23491\end{array}\right ]$ & $43787$\\
\cline{2-9}
	& Lag. & $10^{-2}$ & $[0\,0\,0]$ & $16$ & $48$ & $68$ & $\left [ \begin{array}{c}0\:\:-4.2556\cdot 10^{-1}\:\:-8.9973\cdot 10^{-2}\end{array}\right ]$ & $9.9105\cdot 10^{-1}$\\
	\hline
	\multirow{2}{*}{$\begin{array}{c} \tt UWV \\ n=8 \end{array}$ }& pow. & $1$ &$[0\,0;0\,0]$ & $13$ & $98$ & $188$ &  $\left[\begin{array}{cc}-1.4319\cdot 10^{-5}&-2.6474\cdot 10^{-6}\\ -3.0817\cdot 10^{-1}&-5.6976 \cdot 10^{-2} \end{array}\right]$ & $27.918$\\
\cline{2-9}
    & Lag. & $1$ & $[0\,0;0\,0]$ & $15$ & $65$ & $82$ & $\left[\begin{array}{cc}-1.6755\cdot 10^{-12}&-6.9006\cdot 10^{-13}\\ -3.6060\cdot 10^{-8}&-1.4851 \cdot 10^{-8} \end{array}\right]$ & $1.0000$\\
	\hline
	\multirow{2}{*}{$\begin{array}{c} \tt NN5 \\  n=7\end{array}$ }& pow. & $1$ &$[10\,5]$ & $29$ & $162$ & $300$ & $[12.382\:\:\:9.0331]$ & $3.9610\cdot 10^{-1}$ \\
\cline{2-9}
    & Lag. & $10^{-5}$ & $[10\,5]$ & $18$ & $45$ & $56$ & $[30.931\:\:\:22.295]$ & $1.7652\cdot 10^{-1}$\\
	\hline
	\multirow{2}{*}{$\begin{array}{c} \tt NN1 \\  n=3\end{array}$ }& pow. & $10^{-3}$ &$[0\,30]$ & $15$ & $53$ & $59$ & $[7.9924\:\:\:72.171]$ & $4.2238$ \\
\cline{2-9}
    & Lag. & $10^{-4}$ & $[0\,30]$ & $14$ & $49$ & $52$ & $[26.936\:\:\:177.20]$ & $4.6019$\\
	\hline
	\multirow{2}{*}{$\begin{array}{c} \tt HE1 \\  n=4\end{array}$ }& pow. & $1$ &$[1\,1]$ & $18$ & $73$ & $80$ & $[-1.5482\:\:\:-3.9063]$ & $34.359$ \\
\cline{2-9}
    & Lag. & $10^{-1}$ & $[1\,1]$ & $18$ & $80$ & $87$ & $[-5.1376\:\:\:11.589]$ & $32.168$\\
	\hline
\end{tabular}
\end{center}
\end{table}

In Table \ref{tab:target_poly} we show the influence of the target polynomial on the
computational cost for the $\tt PAS$ system. Open-loop poles
of the system are
\[
\sigma_0 = \left(0,\, 0,\, -9.5970\cdot 10^{-1},\:-36.646\pm523.05i\right)
\]
and we choose alternative target polynomials with the following roots
\begin{align*}
	\sigma_1\:\:=\:\:&\left(-5.0000\cdot 10^{-2},\, -5.0000\cdot 10^{-2},\, -9.5970\cdot 10^{-1},\:-36.646\pm523.05i\right)\\
 	\sigma_2\:\:=\:\:&\left(-1.0000\cdot 10^{-3},\, -1.0000\cdot 10^{-3},\, -9.5970\cdot 10^{-1},\:-36.646\pm523.05i\right)\\
	\sigma_3\:\:=\:\:&\left(0,\, -1.0000\cdot 10^{-4},\, -9.5970\cdot 10^{-1},\:-36.646\pm523.05i\right).
\end{align*}

One can easy to see that the computational cost is decreasing significantly when
the point defining the target polynomial is getting closer to the PENBMI initial iterate.

\begin{table}[h]
\setlength{\tabcolsep}{0.8pt}
\caption {Influence of target polynomial on PENBMI behavior}\label{tab:target_poly}
\begin{center}
\begin{tabular}{|c|c|c|c|c|}
    \hline
    system & \multicolumn{4}{c|}{$\tt PAS$ degree 5} \\
    \hline
    basis&  power & Lagrange & Lagrange & Lagrange\\
    \hline
    $\mu$ & $10^{-3}$ & $10^{-8}$ & $10^{-5}$ & $10^{-2}$\\
 	\hline
    roots & $\sigma_0$ & $\sigma_1$ & $\sigma_2$ & $\sigma_3$\\
    \hline
    $K_0$ & $[0\,0\,0]$ & $[0\,0\,0]$ & $[0\,0\,0]$ & $[0\,0\,0]$\\
    \hline
    out.iter. & $11$ & $19$ & $17$ & $15$\\
    \hline
	inn.iter. & $74$ & $27$ & $33$ & $29$\\
    \hline
	lin.steps & $194$ & $28$ & $44$ & $32$\\
    \hline
	$K^T$ & $\left [ \begin{array}{c}-6.5390\cdot 10^{-4}\\-58.350\\-37.751\end{array}\right ]$ & $\left [ \begin{array}{c}-8.4106\cdot 10^{-6}\\-3.9048\\-9.9675\cdot 10^{-1}\end{array}\right ]$ &    $\left [ \begin{array}{c}-3.3369\cdot 10^{-4}\\-20.480\\-1.2157\end{array}\right ]$ & $\left [ \begin{array}{c}-8.6755\cdot 10^{-8}\\-4.1040\cdot 10^{-1}\\-1.7471\cdot 10^{-1}\end{array}\right ]$\\
    \hline
	$\lambda$ & $73.2917$ & $1.4901\cdot10^{-12}$ & $8.1649\cdot 10^{-3}$ & $2.7241\cdot 10^{-3}$\\
    \hline
    \end{tabular}
\end{center}
\end{table}

Consider the $\tt NN6$ SOF BMI problem that was not solvable in the power basis, see
\cite{cdcecc05}. Open-loop poles of the system are
\[
\sigma_0 = \left(2.7303,\,0,\,-7.2028\cdot 10^{-2}\pm60.804i,\,-1.0785\cdot10^{-1}\pm15.677i,\,-2.6764,\,-3.3000,\,-19.694\right).
\]
The strategy to define the target polynomial is to change the unstable open-loop poles
into slightly stable poles (shifting the real part to a small negative value).
According to this strategy, our target polynomial has the following roots
\begin{align*}
\sigma_1\:\:=\:\:&\left(-1.0000\cdot10^{-3}\pm i,\,-7.2028\cdot 10^{-2}\pm60.804i,\,-1.0785\cdot10^{-1}\pm15.677i,\right.\\
						&\left.\,-2.6764,\,-3.3000,\,-19.694\right).
\end{align*}
The BMI SOF problem is solved with no error or warning in the Lagrange basis,
yielding $\lambda=8.8487\cdot10^{-1}$, $K=[1.3682,\,4.8816,\,44.959,\,59.016]$ with
$17$ outer iterations, $80$ inner iterations and $138$ linesearch steps,
using the orgin as initial point and trade-off parameter $\mu=10^{-5}$.

\section{Conclusion}

The Hermite matrix arising in the symmetric formulation of the polynomial
stability criterion is typically ill-scaled when expressed in the standard
power basis. As a consequence, a nonlinear semidefinite
programming solver such as PENNON may experience convergence problems
when applied on polynomial matrix inequalities (PMIs) coming from benchmark
static output feedback (SOF) problems. In this paper we reformulated Hermite's
SOF PMI in a Lagrange polynomial basis. We slightly extended the results of
\cite{shakoori} to use polynomial interpolation on possibly complex and
repeated nodes to construct the Hermite matrix, bypassing potential numerical
issues connected with Vandermonde matrices. In our control application, a natural
choice of Lagrange nodes are the roots of a target polynomial, the desired
closed-loop characteristic polynomial.

The idea of using the Lagrange polynomial basis to address numerical
problems which are typically ill-scaled when formulated in the power basis
has already proven successful in other contexts. For example, in \cite{fortune}
it was shown that roots of extremely ill-scaled polynomials (such as
a degree 200 Wilkinson polynomial) can be found at machine precision using
eigenvalue computation of generalized companion matrices obtained by
an iterative choice of Lagrange interpolation nodes. In \cite{fft}
the fast Fourier transform (a particular interpolation technique)
was used to perform spectral factorization of polynomials of
degree up to one million. Another example
of successful use of alternative bases and high-degree
polynomial interpolation to address various problems of scientific
computating is the {\tt chebfun} Matlab package, see \cite{chebfun}.
Even though our computational results on SOF PMI problems are less dramatic,
we believe that the use of alternative bases and interpolation can be
instrumental to addressing various other control problems formulated
in a polynomial setting.

\section*{Appendix: Matlab implementation}

A Matlab implementation of the method described in this paper is available
at
\begin{center}
\tt homepages.laas.fr/henrion/software/hermitesof.m
\end{center}
Our implementation uses the Symbolic Math Toolbox and the YALMIP
interface. It is not optimized for efficiency, and therefore
it can be time-consuming already for medium-size examples.

Let us use function {\tt hermitesof} with its default tunings:
\begin{verbatim}
>> [A,B1,B,C1,C] = COMPleib('NN1');
>> A,B,C
A =
     0     1     0
     0     0     1
     0    13     0
B =
     0
     0
     1
C =

     0     5    -1
    -1    -1     0
>> [H,K] = hermitesof(A,B,C)
Quadratic matrix variable 3x3 (symmetric, real, 2 variables)
Linear matrix variable 1x2 (full, real, 2 variables)
\end{verbatim}
Here are some sample entries of the resulting Hermite matrix
\begin{verbatim}
>> sdisplay(H(1,1))
-0.6168744435*K(2)-0.2372594014*K(1)*K(2)+0.04745188027*K(2)^2
>> sdisplay(H(3,2))
0.3019687672*K(1)+0.01984184931*K(2)-0.009656748637*K(1)*K(2)
+0.0003260141644*K(2)^2+0.04013338907*K(1)^2
\end{verbatim}
For this example, the Hermite matrix is quadratic in feedback matrix {\tt K}.
This Hermite matrix is expressed in Lagrange basis, with Lagrange nodes
chosen as the roots of the imaginary part of a random target
polynomial, see the online help of function {\tt hermitesof} for
more information. In particular, it means that each call
to {\tt hermitesof} produces different coefficients. However these
coefficients have comparable magnitudes:
\begin{verbatim}
>> [H,K]=hermitesof(A,B,C);
>> sdisplay(H(1,1))
-0.9592151361*K(2)-0.3689288985*K(1)*K(2)+0.0737857797*K(2)^2
>> sdisplay(H(3,2))
6.702455704*K(1)+0.5150092145*K(2)-0.3440108908*K(1)*K(2)
+0.0184651235*K(2)^2+1.258426367*K(1)^2
\end{verbatim}

The output of function {\tt hermitesof} is reproducible
if the user provides the roots of the target polynomial:
\begin{verbatim}
>> opt = []; opt.roots = [-1 -2 -3];
>> [H,K]=hermitesof(A,B,C,opt);
>> sdisplay(H(1,1))
-0.196969697*K(2)-0.07575757576*K(1)*K(2)+0.01515151515*K(2)^2
>> sdisplay(H(3,2))
0.2*K(1)-0.01818181818*K(2)-0.01212121212*K(1)*K(2)
+0.0007575757576*K(2)^2+0.04166666667*K(1)^2
\end{verbatim}

The Hermite matrix can also be provided in the power basis:
\begin{verbatim}
>> opt = []; opt.basis = 'p';
>> [H,K]=hermitesof(A,B,C,opt);
>> sdisplay(H(1,1))
-13*K(2)-5*K(1)*K(2)+K(2)^2
>> sdisplay(H(3,2))
0
\end{verbatim}

For more complicated examples, the Hermite matrix {\tt H} is not
necessarily quadratic in {\tt K}:
\begin{verbatim}
>> [A,B1,B,C1,C] = COMPleib('NN1');
>> size(B), size(C)
ans =
     5     3
ans =
     3     5
>> [H,K]=hermitesof(A,B,C)
Polynomial matrix variable 5x5 (symmetric, real, 9 variables)
Linear matrix variable 3x3 (full, real, 9 variables)
>> degree(H)
ans =
     5
\end{verbatim}

\section*{Acknowledgments}

This work was partly funded by T\"UBITAK, the Scientific and Technological Research Council of Turkey, and by project No.~103/10/0628 of the Grant Agency of the Czech Republic.

\end{document}